\documentclass[11pt,a4paper]{article}
\usepackage[utf8]{inputenc}
\usepackage[english]{babel}
\usepackage{amsmath,amssymb,amsthm}
\usepackage{geometry}
\usepackage{tcolorbox}
\newtheorem{theorem}{Theorem}[section]
\newtheorem{lemma}[theorem]{Lemma}
\newtheorem{corollary}[theorem]{Corollary}
\theoremstyle{definition}
\newtheorem{definition}[theorem]{Definition}
\newtheorem{remark}[theorem]{Remark}

\numberwithin{equation}{section}

\begin{document}

\title{\textbf{Uniqueness and boundary behaviour of solutions to variational problems with linear growth}} \author{\textbf{Michael Bildhauer} \quad \textbf{Martin Fuchs}} \date{} \maketitle

\begin{abstract}
We investigate the Dirichlet problem for the variational integral $J[u] = \int_{\Omega} f(\nabla u) \, dx$ with density $f$ of linear growth satisfying appropriate ellipticity conditions. We show that the relaxed problem admits a unique solution $u$ in the space of functions of bounded variation, if the set $\Gamma_0$ of convex points $x \in \partial\Omega$ is sufficiently large. For example, the inequality
$\mathcal{H}^{n-1}(\Gamma_0) >
\frac{2}{3}\mathcal{H}^{n-1}(\partial\Omega)$ is sufficient. Moreover, the minimizer $u$ is smooth in the interior of $\Omega$ and attains the prescribed boundary data at least on $\Gamma_0$ in the classical sense.
\end{abstract}


\section{Introduction and statement of the results}

We investigate the minimization problem
\begin{equation}\label{eq:1.1}
J[u] := \int_{\Omega} f(\nabla u) \, dx \to \min , \quad u|_{\partial\Omega} = \varphi, \end{equation} among functions $u: \Omega \rightarrow \mathbb{R}$ with prescribed trace
$\varphi: \partial\Omega \rightarrow \mathbb{R}$ for the case of energy densities $f: \mathbb{R}^n \rightarrow \mathbb{R}$ being of linear growth. If the domain $\Omega$ and also the density are strictly convex and if in addition $\varphi$ is a Lipschitz function, then the Hilbert-Haar approach (see, for instance, \cite{Morrey2008}, Section
4.2) gives the existence of a unique (Lipschitz) solution to problem \eqref{eq:1.1}. In the particular case of nonparametric minimal surfaces, i.e., for the density $f(\nabla u) := \sqrt{1 + |\nabla
u|^2}$, the reader will find a complete discussion concerning existence
and uniqueness of smooth solutions to \eqref{eq:1.1} in the paper \cite{Miranda1971}. For a detailed exposition of the minimal surface case the interested reader is referred to the monograph \cite{Giusti1984}.

Coming back to the general setting, it is a well-known fact that due to the non-reflexivity of the Sobolev space $W^{1,1}(\Omega)$ (see \cite{Adams1975} for a definition of the spaces $W^{k,p}(\Omega)$, $\mathring{W}^{k,p}(\Omega)$, etc.) it does not make sense to replace \eqref{eq:1.1} by \begin{equation}\label{eq:1.2} J[u] \rightarrow \min \quad \text{in } \varphi + \mathring{W}^{1,1}(\Omega).
\end{equation}
This leads towards the relaxation of problems \eqref{eq:1.1} and \eqref{eq:1.2}, respectively, in the space $BV(\Omega)$ of functions with bounded variation (see, e.g., \cite{Giusti1984} for notation). A short discussion of this relaxation procedure will be given below. We emphasize that the BV-approach provides a very elegant existence theory for the relaxed variant of \eqref{eq:1.1} and \eqref{eq:1.2}, respectively, leading to a generalized minimizer having some degree of interior regularity under suitable (ellipticity) assumptions imposed on the density $f$.

However, the questions of uniqueness and of the attainment of the boundary data $\varphi$ remain in general unanswered, and as the example in \cite{Santi1972} shows, positive results can not be expected without further assumptions. This delicate topic is addressed in the work \cite{Beck2018} exhibiting conditions under which actually a unique classical solution of problem \eqref{eq:1.1} exists: suppose that $\Omega$ is a bounded domain of class $C^1$ satisfying an exterior sphere condition and let $\varphi \in C^{1,1}(\bar{\Omega})$. In addition the density $f$ is of radial structure, which means that for $p \in \mathbb{R}^n$ we have $f(p) = g(|p|)$. Then, roughly speaking, the requirement \begin{equation}\label{eq:1.3} \int_1^\infty t \cdot g''(t) \, dt = \infty \end{equation} turns out as a necessary and sufficient condition for the unique solvability of problem \eqref{eq:1.1} in the space $C^{0,1}(\bar{\Omega})$. As a matter of fact, higher interior regularity holds. For a precise formulation we refer to Theorem 1.1 in \cite{Beck2018}. Note that \eqref{eq:1.3} in the context of $\mu$-elliptic integrands (compare inequality \eqref{eq:1.7} below) corresponds to exponents $\mu$ such that \begin{equation}\label{eq:1.4}
1 < \mu < 2.
\end{equation}

In our paper we will first investigate the situation for a larger class of exponents $\mu$, more precisely, condition \eqref{eq:1.4} is replaced by \begin{equation}\label{eq:1.5}
1 < \mu < \begin{cases} 3, & \text{if } n \geq 3 \\ \leq 3, & \text{if } n = 2 \end{cases}, \end{equation} and as already remarked before, neither uniqueness nor continuous attainment of the boundary data can be expected. A second topic of our note concerns energy densities $f$ of linear growth with the particular structure \[ f(\nabla u) = g_1(|\partial_1 u|) + \dots + g_n(|\partial_n u|), \] which do not fall in the category of \mbox{$\mu$-elliptic} integrands.
In both cases we will exhibit geometric conditions implying uniqueness and continuity up to the boundary. Before stating these results we fix our

\subsection*{Notation and assumptions:}
Let the following general assumptions (GA) hold:
\begin{itemize}
 \item[\textbf{(GA1)}] $\Omega \subset \mathbb{R}^n$, $n \geq 2$, is a bounded domain with Lipschitz continuous boundary.
 \item[\textbf{(GA2)}] The boundary datum $\varphi$ is in the space $W^{1,2}(\Omega)$ with Hölder continuous trace on $\partial\Omega$.
 \item[\textbf{(GA3)}] The energy density $f: \mathbb{R}^n \rightarrow \mathbb{R}$ is of class $C^2(\mathbb{R}^n)$ with second partial derivatives being locally Hölder continuous. We assume that $f$ is of linear growth in the sense that
 \begin{equation}\label{eq:1.6}
 a(|p| - 1) \leq f(p) \leq A(|p| + 1), \quad p \in \mathbb{R}^n,
 \end{equation}
 holds with constants $a, A > 0$.

 The density $f$ is $\mu$-elliptic, which means that there exist numbers $\mu > 1$, $c_1, c_2 > 0$ such that
 \begin{equation}\label{eq:1.7}
 c_1(1 + |p|^2)^{-\mu/2} |q|^2 \leq D^2f(p)(q,q) \leq c_2(1 +
|p|^2)^{-1/2} |q|^2
 \end{equation}
 holds for all $p, q \in \mathbb{R}^n$.
\end{itemize}
Note that in the case $f(p) = \sqrt{1 + |p|^2}$ of nonparametric minimal surfaces we have \eqref{eq:1.6} and \eqref{eq:1.7} with the choice \mbox{$\mu = 3$}.

Let us summarize some general facts concerning the relaxation of problem \eqref{eq:1.1} in the space $BV(\Omega)$. The interested reader is referred for instance to \cite{Giusti1984, Beck2018, Ambrosio2000, Giaquinta1979} and the references quoted therein. We start by introducing the quadratic regularization of problem \eqref{eq:1.1}: for
$0 < \delta < 1$ consider the problem
\begin{equation}\label{eq:1.8}
\begin{cases}
J_\delta[v] := \frac{\delta}{2} \int_{\Omega} |\nabla v|^2 \, dx + J[v] \\ \rightarrow \min \text{ in } \varphi + \mathring{W}^{1,2}(\Omega), \quad J[v] = \int_{\Omega} f(\nabla v) \, dx, \end{cases} \end{equation} and recall the following facts:

\begin{lemma}
Under the assumptions (GA1-3) we have
\begin{itemize}
 \item[i)] Problem \eqref{eq:1.8} admits a unique solution $u_\delta$.
 \item[ii)] The minimizer is of class $W_{\text{loc}}^{2,2}(\Omega) \cap C^1(\Omega)$ satisfying
 \begin{equation}\label{eq:1.9}
 \|u_\delta\|_{L^\infty(\Omega)} \leq \max_{\partial\Omega}
|\varphi|.
 \end{equation}
 \item[iii)] The regularity of $u_\delta$ can be improved up to $u_\delta \in C^2(\Omega) \cap C^{0,\alpha}(\bar{\Omega})$ for some $\alpha \in (0,1)$.
 \item[iv)] (comparison principle) Consider an affine function $\ell:
\mathbb{R}^n \rightarrow \mathbb{R}$ such that $\varphi \leq \ell$ on $\partial\Omega$ (or $\geq$). Then it holds
 \begin{equation}\label{eq:1.10}
 u_\delta(x) \leq \ell(x) \quad \text{for all } x \in \bar{\Omega} \ (\text{or } \geq).
 \end{equation}
\end{itemize}
\end{lemma}

The assertions i) and ii) are immediate, for $u_\delta \in C^2(\Omega)$ we refer to interior Campanato type estimates (using the local Hölder continuity of $D^2f$), and $u_\delta \in C^{0,\alpha}(\bar{\Omega})$ follows for instance from Theorem 7.8, p.232, in \cite{Giusti2003}.
Finally, statement iv) is a consequence of Theorem 10.1, p.263, from \cite{Gilbarg2001}.

Next we summarize the results of Theorem 4.14 and 4.16 in \cite{Bildhauer2003} observing that \eqref{eq:1.9} implies Assumption
4.11 required in \cite{Bildhauer2003}.

\begin{theorem}
Let (GA1-3) hold and assume that we have \eqref{eq:1.5} for the exponent $\mu$ introduced in \eqref{eq:1.7}. We let \[
\mathcal{M}_0 := \{u^* \in BV(\Omega) : u^* \text{ is a } L^1(\Omega)\text{-cluster point of the family } (u_\delta)_{0<\delta<1} \text{ introduced in Lemma 1.1}\}, \] \[ \mathcal{M} := \{u^* \in BV(\Omega) : u^* \text{ is a } L^1(\Omega)\text{-cluster point of a } J\text{-minimizing sequence in } \varphi + \mathring{W}^{1,1}(\Omega)\}.
\]
Then it holds:
\begin{itemize}
 \item[i)] $\emptyset \neq \mathcal{M}_0 \subset \mathcal{M}$, and any $u \in \mathcal{M}$ minimizes the functional
 \begin{equation}\label{eq:1.11}
 K[v] := \int_{\Omega} f(\nabla v) + \int_{\partial\Omega} f_\infty((\varphi - v)\mathcal{N}) \, d\mathcal{H}^{n-1}
 \end{equation}
 in the class $BV(\Omega)$, moreover, it holds
 \[
 \inf_{BV(\Omega)} K = \inf_{\varphi + \mathring{W}^{1,1}(\Omega)} J.
 \]
\end{itemize}
In \eqref{eq:1.11} the quantity $f(\nabla v)$ has to be understood as a convex function of a measure. $f_\infty$ is the recession function of $f$, $\mathcal{N}$ denotes the exterior normal of $\partial\Omega$, and in $\int_{\partial\Omega} \dots$ the symbol $v$ stands for the BV-trace of $v$ on $\partial\Omega$. Let us fix some function $u^* \in \mathcal{M}_0$. We have
\begin{itemize}
\item[ii)] $u^* \in C^2(\Omega) \cap W^{1,1}(\Omega)$, if $\mu < 3$. In the case $n=2$ together with $\mu=3$ we just have
$u^* \in W^{1,1}(\Omega)$ together with some local higher integrability of $\nabla u^*$.
\item[iii)] For any $u \in \mathcal{M}$ it holds $\nabla u = \nabla u^*$ yielding uniqueness up to additive constants.
\item[iv)] The comparison principle iv) from Lemma 1.1 is valid for $u^*$.
\end{itemize}
\end{theorem}
Concerning iv) we note that this claim immediately follows from \eqref{eq:1.10} and the choice of $u^*$. Unfortunately Theorem 1.2 says nothing about the attainment of the boundary values $\varphi$ and we just have uniqueness up to constants. In order to get some insight we introduce the following subsets of $\partial\Omega$:\begin{definition}Let the domain $\Omega$ satisfy (GA1) and assume that $\varphi$ is just a function in $C^0(\partial\Omega)$. We then define \\\\\(\Gamma _{0}:=\Gamma _{0}(\partial \Omega ):=\{x_{0}\in \partial \Omega :\text{there\ exists\ an\ open\ ball\ }B\text{\ such\ that\ }\Omega \subset B\text{\ and\ }x_{0}\in \partial B\},\)\
\begin{align*}\Gamma := \Gamma(\partial\Omega, \varphi)
:= \{x_0 \in \partial\Omega :
&\text{ for any } \varepsilon > 0 \text{ we find } a,b\in\mathbb{R}^n \
&\text{such that for all } x\in\partial\Omega \
&\text{ it holds }\\
\varphi(x_0)-\varepsilon+a\cdot(x-x_0)
\le \varphi(x) 
&\le \varphi(x_0)+\varepsilon+b\cdot(x-x_0)
\}.
\end{align*}
\end{definition}
The points $x_0 \in \partial\Omega$ being contained in $\Gamma_0$ can be seen as points of convexity of $\partial\Omega$
to be understood in a rather strong sense. It holds
\begin{lemma}
\label{lem:1.4}For $\Omega$ and $\varphi$ as in Definition 1.3 we have $\emptyset \neq \Gamma_0 \subset \Gamma$.
\end{lemma}
Using Lemma \ref{lem:1.4} we will show
\begin{theorem}\label{thm:1.5}Let (GA1-3) hold together with condition \eqref{eq:1.5} for the exponent $\mu$. Consider $u^* \in \mathcal{M}_0$ and a point $x_0 \in \Gamma$. We
have
\begin{itemize}
\item[i)] $\displaystyle \lim_{\Omega \ni x \rightarrow x_0} u^*(x) = \varphi(x_0)$, if $\mu < 3$,
\item[ii)] $\displaystyle \operatorname*{ess \, sup}_{\Omega \cap B\rho(x_0)} |u^*(x) - \varphi(x_0)| \rightarrow 0$ as $\rho \rightarrow 0$, if the case $\mu = 3$ together with $n = 2$ is considered.
\end{itemize}
\end{theorem}
\begin{corollary}\label{cor:1.6}
Let
(GA1-3) together with \eqref{eq:1.5} hold. Then $\mathcal{M}_0$ contains exactly one element $u^*$. On the boundary part $\Gamma$ the BV-trace of $u^*$ equals $\varphi$ $\mathcal{H}^{n-1}$-almost everywhere, more precisely, i) and ii) from Theorem \ref{thm:1.5} hold for points $x_0 \in \Gamma$.\end{corollary}If a generalized minimizer $u \in \mathcal{M}$ takes the boundary value $\varphi(x_0)$ at some point $x_0 \in \Gamma$ in the sense that i) or ii) from Theorem \ref{thm:1.5} hold for the function $u$, then clearly (recall Theorem 1.2 iii)) we must have that $u = u^*$ with $u^*$ from Corollary \ref{cor:1.6}. Thus the attainment of the boundary data at least for one point $x_0 \in \Gamma$ guarantees uniqueness in the sense that $\mathcal{M}_0 = \mathcal{M} = \{u^*\}$.
However, there still might exist minimizers $\tilde{u} \in \mathcal{M}$ of the form $\tilde{u} = u^* + c$ with $u^*$ from Corollary \ref{cor:1.6} and for a number $c \neq 0$. This situation can be excluded under appropriate assumptions on the size of $\Gamma$.\begin{theorem}\label{thm:1.7}Let (GA1-3) hold together with \eqref{eq:1.5} and assume in addition that the density $f$ has a radial structure, which means that $f(p) = g(|p|)$ for a suitable function $g$.
Assume in addition that the set $\Gamma$ introduced in Definition 1.3
satisfies\begin{equation}\label{eq:1.12}\mathcal{H}^{n-1}(\Gamma) > \frac{2}{3} \mathcal{H}^{n-1}(\partial\Omega).\end{equation}Then the set $\mathcal{M}$ of $K$-minimizers (compare \eqref{eq:1.11}) contains exactly one element $u$ for which the boundary condition in the sense of
i) and ii) of Theorem \ref{thm:1.5} is satisfied in any point $x_0 \in
\Gamma$.\end{theorem}
\begin{remark}i) We can replace \eqref{eq:1.12} by the stronger but purely geometric condition (not involving
$\varphi$)
\begin{equation}\label{eq:1.13}
\mathcal{H}^{n-1}(\Gamma_0) > \frac{2}{3} \mathcal{H}^{n-1}(\partial\Omega).
\end{equation}
For instance, if we let $n = 2$ and take $\Omega := B_R(x) - \overline{B_r(x)}$ with radii $0 < r < R$, then it holds
$\Gamma_0 = \partial B_R(x)$ and \eqref{eq:1.13} just means that $R > 2r$.
ii) Returning to problem \eqref{eq:1.1} and \eqref{eq:1.2}, respectively, we have shown that under the assumptions of Theorem \ref{thm:1.7} a unique function $u \in W^{1,1}(\Omega) (\cap C^2(\Omega))$ can be found for which $u = \varphi$ holds on $\Gamma$ and which is $J$-minimizing in an appropriate sense.
\end{remark}Next we consider densities without radial structure, for example we let\begin{equation}\label{eq:1.14}f(p) =
g_1(|p_1|) + \dots + g_n(|p_n|), \quad p \in \mathbb{R}^n,\end{equation}with functions $g_i: [0,\infty) \rightarrow \mathbb{R}$ being of class $C^{2,\alpha}([0, T])$ for each $T < \infty$ with some $\alpha$ eventually depending on $T$ and in addition
satisfying
\begin{equation}\label{eq:1.15}
\begin{cases}
g_i'(0) = 0, \quad g_i' \in L^\infty([0,\infty)) \text{ and} \\
c_i(1+t)^{-\mu_i} \leq g_i''(t) \leq C_i
\text{ for all } t \geq 0
\end{cases}
\end{equation}
with constants $c_i, C_i > 0$ and exponents $\mu_i > 1$. Note that \eqref{eq:1.15} implies \eqref{eq:1.6} in (GA3) for $f$ defined in \eqref{eq:1.14}, but in place of \eqref{eq:1.7} we just get \(c(1+|p|^{2})^{-\mu /2}|q|^{2}\le D^{2}f(p)(q,q)\le C|q|^{2}\) for $p, q \in \mathbb{R}^n$, with positive constants $c, C$ and for the choice $\mu := \max\{\mu_i : i = 1, \dots, n\}$. The setting described in \eqref{eq:1.14} and \eqref{eq:1.15} is covered by\begin{theorem}Let
(GA1-3) hold, where in (GA3) inequality \eqref{eq:1.7} is replaced by\begin{equation}\label{eq:1.16}\lambda(1+|p|^2)^{-\mu/2} |q|^2 \leq
D^2f(p)(q,q) \leq \Lambda(1+|p|^2)^{-\varkappa/2} |q|^2, \quad p, q \in \mathbb{R}^n,\end{equation}for constants $\lambda, \Lambda > 0$ and with exponents $\mu$ and $\varkappa$ such that\begin{equation}\label{eq:1.17}\mu > 1, \quad 0 \leq \varkappa \leq 1, \quad \mu < 2 + \varkappa.\end{equation}Moreover, we assume that $f \geq 0$ and $f(-p) = f(p)$, $p \in \mathbb{R}^n$, hold. Then we
have:\begin{itemize}\item[i)] The results of Theorem 1.2, Theorem \ref{thm:1.5} and of Corollary \ref{cor:1.6} continue to hold, in particular we have $u^* \in C^2(\Omega) \cap W^{1,1}(\Omega)$ with $u^*$ from Theorem 1.2.\item[ii)] Consider the measure ($A \subset \partial\Omega$) \\\\\(m(A):=\int _{A}f_{\infty }(\mathcal{N})\,d\mathcal{H}^{n-1},\) $\\\\\mathcal{N}$ denoting the exterior normal of $\partial\Omega$ and $f_\infty$ being the recession function of $f$. Suppose that we know\begin{equation}\label{eq:1.18}m(\Gamma) > \frac{2}{3} m(\partial\Omega)\end{equation}for the set $\Gamma$ introduced in Definition 1.3. Then the set $\mathcal{M}$ introduced in Theorem 1.2 consists of exactly one element $u$, which means that the relaxed problem \[K\rightarrow \min \quad \text{in\ }BV(\Omega )\] with $K$ from \eqref{eq:1.11} has a unique solution, which in addition satisfies i) from Theorem
\ref{thm:1.5}.\end{itemize}\end{theorem}\begin{remark}i) If we consider functions $g_i \geq 0$ as in \eqref{eq:1.15} and define $f$ according to \eqref{eq:1.14}, then \eqref{eq:1.16} holds with the choices $\varkappa = 0$, $\mu = \max\{\mu_1, \dots, \mu_n\}$, thus \eqref{eq:1.17} is satisfied for $1 < \mu_i < 2$. The measure $m$ takes the form ($A \subset \partial\Omega$) \\\\\(m(A)=\int _{A}\sum _{i=1}^{n}|\mathcal{N}_{i}|\lim _{t\rightarrow \infty }g_{i}^{\prime
}(t)\,d\mathcal{H}^{n-1}.\)\\\\
ii) Clearly we can replace \eqref{eq:1.18} by the stronger condition $m(\Gamma_0) > \frac{2}{3} m(\partial\Omega)$.\end{remark}\section{Some preliminary geometric remarks -- proof of Lemma \ref{lem:1.4}}Let $\Omega$ and $\varphi$ satisfy the assumptions of Lemma \ref{lem:1.4} and choose a ball $B_R(\eta)$ such that $\bar{\Omega} \subset B_R(\eta)$. If $R^*$ denotes the infimum of the radii of balls centered at $\eta$ having this property,
then it holds
\begin{equation}\label{eq:2.1}
\Omega \subset B_{R^*}(\eta).
\end{equation}
If this is not the case, we find $x \in {\Omega}$ with the property $|x - \eta| \geq R^*$. Since clearly $\Omega \subset \overline{B_{R^*}(\eta)}$, it holds $|x - \eta| = R^*$. The set $\Omega$ is open, thus $B_r(x) \subset \Omega$ for some radius $r > 0$, and in conclusion $B_r(x) \subset \overline{B_{R^*}(\eta)}$ leading to a contradiction, hence \eqref{eq:2.1} is true. Next we claim
\begin{equation}\label{eq:2.2}\partial\Omega \cap \partial
B_{R^*}(\eta) \neq \emptyset.
\end{equation}If \eqref{eq:2.2} would be wrong, we get from \eqref{eq:2.1} that $\bar{\Omega} \subset B_{R^*}(\eta)$, hence $\bar{\Omega} \subset B_{R^* -\varepsilon}(\eta)$ for $\varepsilon > 0$ sufficiently small, contradicting the minimality of $R^*$. From \eqref{eq:2.1} and \eqref{eq:2.2} we deduce that $\Gamma_0 \neq \emptyset$. In the next step we prove the inclusion $\Gamma_0 \subset \Gamma$. Let $x_0 \in \Gamma_0$ and choose a ball $B_R(y_0)$ such that\begin{equation}\label{eq:2.3}\Omega \subset B_R(y_0), \quad
x_0 \in \partial\Omega \cap \partial B_R(y_0).\end{equation}For $t > 0$ sufficiently large we let \(\eta _{t}:=x_{0}+t(y_{0}-x_{0}),\quad r_{t}:=t|y_{0}-x_{0}|=tR\) and consider the ball $B_{r_t}(\eta_t)$ containing $B_R(y_0)$ and satisfying (recall \eqref{eq:2.3}) $\partial
B_{r_t}(\eta_t) \cap \partial B_R(y_0) = \{x_0\}$. We also observe that\begin{equation}\label{eq:2.4}(\eta_t - x_0) \cdot (y - x_0) = t(y_0
- x_0) \cdot (y - x_0) \geq 0\end{equation}for points $y \in \overline{B_R(y_0)}$. Given $\varepsilon > 0$ the continuity of $\varphi$ at $x_0$ gives a number $\delta > 0$ such
that
\begin{equation}\label{eq:2.5}\varphi(x) \leq \varphi(x_0) + \varepsilon, \quad x \in \partial\Omega \cap B_\delta(x_0),\end{equation}and from \eqref{eq:2.4}, \eqref{eq:2.5} it follows (observe $\bar{\Omega} \subset
\overline{B_R(y_0)}$)\begin{equation}\label{eq:2.6}\varphi(x) \leq
\varphi(x_0) + \varepsilon + (\eta_t - x_0) \cdot (x - x_0)\end{equation}for any $x \in \partial\Omega \cap B_\delta(x_0)$. For points $x \in \overline{B_R(y_0) - B_\delta(x_0)}$ we have \[(x-x_{0})\cdot (\eta _{t}-x_{0})=(x-x_{0})\cdot (y_{0}-x_{0})t\geq c\,t\,|x-x_{0}||y_{0}-x_{0}|\geq c\,t\,\delta \,R,\] where $c$ is a positive constant independent of $x$. Thus, by enlarging $t$, if necessary, we can achieve that by letting $b := \eta_t - x_0$ it holds \[b\cdot (x-x_{0})\ge \max_{\partial \Omega }\varphi -\varphi (x_{0})-\varepsilon \] for all $x \in \overline{B_R(y_0) - B_\delta(x_0)}$, in particular we
get\begin{equation}\label{eq:2.7}\varphi(x) \leq \varphi(x_0) + \varepsilon + b \cdot (x - x_0)\end{equation}for points $x \in \partial\Omega$, $|x - x_0| \geq \delta$. From \eqref{eq:2.6} and \eqref{eq:2.7} we finally get the upper bound $\varphi(x) \leq
\varphi(x_0) + \varepsilon + b \cdot (x - x_0)$, $x \in \partial\Omega$, with vector $b$ from above. The existence of $a \in \mathbb{R}^n$ with the property that $\varphi(x) \geq \varphi(x_0) - \varepsilon + a \cdot (x - x_0)$ holds for all $x \in \partial\Omega$ follows in an analogous manner, which completes the proof of Lemma \ref{lem:1.4}, since we have shown that $x_0 \in \Gamma$ is an element of $\Gamma$. \qed
\section{Proofs of Theorem 1.5, Corollary 1.6 and Theorem 1.7}We start with the\subsection*{Proof of Theorem \ref{thm:1.5}.}Let the assumptions of this theorem hold and choose a function $u^* \in \mathcal{M}_0$, which by the definition of the set $\mathcal{M}_0$ stated in Theorem 1.2 means that $u^* = \lim_{\delta\rightarrow0} u_\delta$ with functions $u_\delta$ introduced in Lemma 1.1 and where ``$\lim_{\delta\rightarrow0}$'' has to be understood as the $L^1$-limit of some sequence $(u_{\delta_k})$ with $\delta_k \rightarrow 0$. Fix some point $x_0 \in \Gamma$ and let $\varepsilon > 0$ be given. According to Definition 1.3 there exist $a, b \in \mathbb{R}^n$ with the property that\begin{equation}\label{eq:3.1}a \cdot (x - x_0) - \varepsilon +
\varphi(x_0) \leq \varphi(x) \leq b \cdot (x - x_0) + \varepsilon + \varphi(x_0)\end{equation}holds for all $x \in \partial\Omega$. Quoting inequality \eqref{eq:1.10} from Lemma 1.1 iv) we deduce from \eqref{eq:3.1}\begin{equation}\label{eq:3.2}a \cdot (x - x_0) - \varepsilon \leq u_\delta(x) - \varphi(x_0) \leq b \cdot (x - x_0) + \varepsilon\end{equation}valid for $x \in \bar{\Omega}$. Since we may assume that $u_\delta(x) \rightarrow u^*(x)$ holds for almost all $x \in \Omega$, \eqref{eq:3.2} yields\begin{equation}\label{eq:3.3}a \cdot (x -
x_0) - \varepsilon \leq u^*(x) - \varphi(x_0) \leq b \cdot (x - x_0) + \varepsilon\end{equation}for almost every point $x \in \Omega$. If the case $\mu < 3$ is considered, then we can quote Theorem 1.2 ii), hence \eqref{eq:3.3} is true for any $x \in \Omega$ and \eqref{eq:3.3} yields part i) of Theorem \ref{thm:1.5}. For the limit case $\mu = 3$ together with $n = 2$ we note that $u^*$ is in $W^{1,1}(\Omega)$ (see again Theorem 1.2 ii)), so that now \eqref{eq:3.3} implies ii) of Theorem \ref{thm:1.5}. \qed
\subsection*{Proof of Corollary \ref{cor:1.6}.}For $u^*, \tilde{u} \in \mathcal{M}_0$ it follows from Theorem 1.2 i) and
iii) that $\tilde{u} = u^* + c$ for some number $c \in \mathbb{R}$.
According to Lemma \ref{lem:1.4} we can select a point $x_0 \in \Gamma_0 \subset \Gamma$, and on account of Theorem \ref{thm:1.5} we see that only $c = 0$ is possible. In conclusion, the set $\mathcal{M}_0$ contains exactly one element $u^*$ and for any sequence $\delta_k \rightarrow 0$ it holds $\lim_{k\rightarrow\infty} u_{\delta_k} = u^*$ in $L^1(\Omega)$.
Quoting formula (2.13) from Theorem 2.10 in \cite{Giusti1984} we find on account of the description of the limit behaviour of $u^*$ in points $x_0 \in \Gamma$ given in Theorem \ref{thm:1.5} that the trace of $u^*$ equals $\varphi$ in $\mathcal{H}^{n-1}$-almost all points of $\Gamma$.
\qed\subsection*{Proof of Theorem \ref{thm:1.7}.}Let $u^*$ denote the unique generalized minimizer in the class $\mathcal{M}_0$ (see Corollary
\ref{cor:1.6}) and consider $u \in \mathcal{M}$. From Theorem 1.2 iii) it follows $u = u^* + c$ for some $c \in \mathbb{R}$. Let us assume that $c > 0$. Then it holds due to $u^* = \varphi$ $\mathcal{H}^{n-1}$-almost everywhere on
$\Gamma$\begin{equation}\label{eq:3.4}\mathcal{H}^{n-1}(\Gamma) \leq \mathcal{H}^{n-1}(\{x \in \partial\Omega : u(x) > \varphi(x)\}).\end{equation}Quoting \cite{Bildhauer2018}, Corollary 2.5, we find\begin{equation}\label{eq:3.5}\mathcal{H}^{n-1}([u > \varphi]) \leq \mathcal{H}^{n-1}([u = \varphi]) + \mathcal{H}^{n-1}([u < \varphi]),\end{equation}where $[u > \varphi] := \{x \in \partial\Omega :
u(x) > \varphi(x)\}$, etc., and the symbol $u$ also denotes the trace on $\partial\Omega$. Observing the inclusions (valid up to sets of $\mathcal{H}^{n-1}$-measure zero) \([u=\varphi ]\subset \partial \Omega -\Gamma ,\quad [u<\varphi ]\subset \partial \Omega -\Gamma ,\) which again follow from $u^* = \varphi$ $\mathcal{H}^{n-1}$-almost everywhere together with $c > 0$, we get from \eqref{eq:3.5}\begin{equation}\label{eq:3.6}\mathcal{H}^{n-1}([u >
\varphi]) \leq 2 \left( \mathcal{H}^{n-1}(\partial\Omega) -
\mathcal{H}^{n-1}(\Gamma) \right).\end{equation}With \eqref{eq:3.6} we return to \eqref{eq:3.4} with the result \(\mathcal{H}^{n-1}(\Gamma )\le 2\mathcal{H}^{n-1}(\partial \Omega )-2\mathcal{H}^{n-1}(\Gamma ).\) So if we assume the validity of inequality \eqref{eq:1.12}, then the equation $u = u^* + c$ with some positive constant $c$ leads to a contradiction.
The case $c < 0$ is discussed in the same manner yielding the claim of Theorem \ref{thm:1.7}. Note that the radial structure of the density $f$ just enters through assumption (7) made in \cite{Bildhauer2018} leading to inequality \eqref{eq:3.5} from above. \qed\section{Proof of Theorem 1.9}Let the assumptions of Theorem 1.9 hold. Then the results collected in i) are consequences of Corollary 1.2 in \cite{Bildhauer2020}, where even the case $\varkappa \in (-1,0]$ is considered. So we focus on the proof of ii). As in the proof of Theorem \ref{thm:1.7} we denote by $u^*$ the unique element of $\mathcal{M}_0$ and consider some $u \in \mathcal{M}$, which must take the form $u = u^* + c$. In order to prove that $c = 0$, we recall (see Theorem 1.2 i)) that $K[u] \leq K[u+t]$ holds for all $t \in \mathbb{R}$.The definition of the functional $K$ (see \eqref{eq:1.11}) shows that the function (recall $f \geq 0$ and
$f(-p) = f(p)$) \(g:\mathbb{R}\rightarrow [0,\infty ),\quad g(t):=\int _{\partial \Omega }f_{\infty }((u+t-\varphi )\mathcal{N})\,d\mathcal{H}^{n-1},\) attains its minimum at $t = 0$. We
have\begin{align*}
g(t) &= \int_{\partial\Omega} |u+t-\varphi|
f_\infty(\mathcal{N}) \, d\mathcal{H}^{n-1} \\
&= |t| m([u = \varphi]) + \int_{\partial\Omega \cap [u \neq \varphi]} |u+t-\varphi| \, dm =: g_0(t)
+ g_1(t),
\end{align*}where as usual we have set $[u = \varphi] := \{x \in
\partial\Omega : u(x) = \varphi(x)\}$, etc., and $m$ denotes the measure
$m(A) := \int_A f_\infty(\mathcal{N}) \, d\mathcal{H}^{n-1}$, $A \subset \partial\Omega$. The decomposition $g = g_0 + g_1$ together with $0 \in \partial g(0)$ implies according to (10) in
\cite{Bildhauer2018}\begin{equation}\label{eq:4.1}0 \in [\eta - \alpha, \eta + \alpha]\end{equation}with \mbox{$\alpha := m([u = \varphi])$}, \mbox{$\eta := \int_{\partial\Omega \cap [u \neq \varphi]} \frac{u-\varphi}{|u-\varphi|} \, dm$}. Clearly \eqref{eq:4.1} can be rewritten in the form\begin{equation}\label{eq:4.2}|m([u > \varphi]) - m([u < \varphi])| \leq m([u = \varphi]),\end{equation}which corresponds to inequality \eqref{eq:3.5}. Starting from \eqref{eq:4.2}, we can repeat the proof of Theorem \ref{thm:1.7} replacing $\mathcal{H}^{n-1}$ by the measure $m$ ending up with $c = 0$, hence $u = u^*$, if we use condition \eqref{eq:1.18}.
\qed

\vfill

\begin{minipage}[t]{0.48\textwidth}

\textbf{Michael Bildhauer}\\
Saarland University\\
Department of Mathematics\\
P.O. Box 15 11 50\\
66041 Saarbrücken, Germany\\
\texttt{bibi@math.uni-sb.de}

\end{minipage}
\hfill
\begin{minipage}[t]{0.48\textwidth}

\textbf{Martin Fuchs}\\
Saarland University\\
Department of Mathematics\\
P.O. Box 15 11 50\\
66041 Saarbrücken, Germany\\
\texttt{fuchs@math.uni-sb.de}

\end{minipage}

\newpage
\begin{thebibliography}{99}\bibitem{Morrey2008}
C.B. Morrey, \textit{Multiple integrals in the calculus of variations}.
Springer Verlag, Berlin 2008.\bibitem{Miranda1971} M. Miranda, Un principio di massimo forte per le frontiere minimali e una sua applicazione alla risoluzione del problema al contorno per l'equazione delle superfici di area minima. \textit{Rend. Sem. Mat. Univ. Padova} 45 (1971), 355-366.\bibitem{Giusti1984} E. Giusti, \textit{Minimal surfaces and functions of bounded variation}. Birkhäuser, Basel 1984.\bibitem{Adams1975} R.A. Adams, \textit{Sobolev spaces}, Academic Press, New York 1975.\bibitem{Santi1972} E. Santi, Sul problema al contorno per l'equazione delle superfici di area minima su domini limitati qualunque. \textit{Ann. Univ. Ferrara, N. Ser., Sez. VII} 17 (1972), 13-26.\bibitem{Beck2018} L. Beck, M. Bulíček, E. Maringová, Globally Lipschitz minimizers for variational problems with linear growth. \textit{ESAIM Control Optim. Calc. Var.} 24 (2018), no. 4, 1395-1413.\bibitem{Ambrosio2000} L. Ambrosio, N. Fusco, D. Pallara, \textit{Functions of bounded variation and free discontinuity problems}.
Oxford University Press, New York 2000.\bibitem{Giaquinta1979} M.
Giaquinta, G. Modica, J. Souček, Functionals with linear growth in the calculus of variations. I. \textit{Comm. Math. Univ. Carolinae} 20 (1979), no. 1, 143-156.\bibitem{Giusti2003} E. Giusti, \textit{Direct methods in the calculus of variations}. World Scientific Publishing, Singapore 2003.\bibitem{Gilbarg2001} D. Gilbarg, N.S. Trudinger, \textit{Elliptic partial differential equations of second order}, Springer Verlag, Berlin 2001.\bibitem{Bildhauer2003} M. Bildhauer, \textit{Convex variational problems}, Lecture Notes in Math. 1818, Springer Verlag, Berlin 2003.\bibitem{Bildhauer2018} M. Bildhauer, M.
Fuchs, Some remarks on the (non-) attainment of the boundary data for variational problems in the space BV. \textit{J. Convex Anal.} 25 (2018), 219-223.\bibitem{Bildhauer2020} M. Bildhauer, M. Fuchs, Splitting type variational problems with linear growth conditions, \textit{J. Math. Sci. (N.Y.)} 250, no. 2 (2020), 232-249.\end{thebibliography}
\end{document}